\documentclass[11pt]{article}
\usepackage[russian, english]{babel}
\usepackage[cp1251]{inputenc}
\usepackage{amssymb,amsmath}
\def\abd{\mathop {\rm \leftharpoondown \!\! \rightharpoondown} }
\begin{document}
\begin{center}
\textbf{THE DISCRETE ANALOGUE OF THE OPERATOR
$\displaystyle\frac{d^{2m}}{dx^{2m}}$ AND \\
 ITS PROPERTIES}\\
\textbf{ Kh.M.Shadimetov}
\end{center}

\begin{abstract}
In this paper the discrete analogue $D_m[\beta]$ of the
differential operator $d^{2m}/dx^{2m}$ is constructed and its some
new properties are proved.
\end{abstract}


\section{Main results.}

First S.L.Sobolev [1] studied construction and investigated
properties of the operator $ D_{hH}^{(m)}[\beta]$, which is
inverse of the convolution operator with function
$G_{hH}^{(m)}[\beta]=h^nG_m(hH\beta)$. The function
$D_{hH}^{(m)}[\beta]$ of discrete variable, satisfying the
equality
$$
h^n D_{hH}^{(m)}[\beta]*G_{hH}^{(m)}[\beta]=\delta[\beta]
$$
is called by \emph{the discrete analogue} of the polyharmonic
operator $\Delta^m$. S.L.Sobolev suggested an algorithm for
finding function $D_{hH}^{(m)}[\beta]$ and proved several
properties of this function. In one dimensional case, i.e. the
discrete analogue of the operator
$\displaystyle\frac{d^{2m}}{dx^{2m}}$ was constructed by
Z.Zh.Zhamalov [2, 3]. But there the form of this function was
written with $m+1$ unknown coefficients. In works [4, 5] these
coefficients were found, hereunder the discrete analogue of the
operator $\displaystyle\frac{d^{2m}}{dx^{2m}}$ was constructed
completely.

In this paper we give the results of works  [4-7], concerning to
construction of the discrete analogue $D_m[\beta]$ of the operator
$\displaystyle\frac{d^{2m}}{dx^{2m}}$, and discovery of its
properties, which early were not known.

Following statements are valid.

{\bf Theorem 1.} {\it The discrete analogue of the differential
operator $\displaystyle\frac{d^{2m}}{dx^{2m}}$ have following form
$$
D_m[\beta]=\frac{(2m-1)!}{h^{2m}}\left\{
\begin{array}{lll}
{\displaystyle
\sum\limits_{k=1}^{m-1}\frac{(1-\lambda_k)^{2m+1}\lambda_k^{|\beta|}}
{\lambda_kE_{2m-1}(\lambda_k)}        }
& \mbox{ for }& |\beta|\geq 2,\\
{\displaystyle
1+\sum\limits_{k=1}^{m-1}\frac{(1-\lambda_k)^{2m+1}}
{E_{2m-1}(\lambda_k)}                  }
& \mbox{ for }& |\beta|= 1,\\
{\displaystyle
-2^{2m-1}+\sum\limits_{k=1}^{m-1}\frac{(1-\lambda_k)^{2m+1}}
{\lambda_k E_{2m-1}(\lambda_k)}           } & \mbox{ for }& \beta=
0,
\end{array}
\right. \eqno (1)
$$
where $E_{\alpha}(\lambda)$ is the Euler polynomial of degree
$\alpha$, $\lambda_k$ are the roots of the Euler polynomial
$E_{2m-2}(\lambda)$, in module less than unity, i.e.
$|\lambda_k|<1$,  $h$ is the step of the lattice.}

{\bf Property 1.} {\it  The discrete analogue $D_m[\beta]$ of the
differential operator of order $2m$ have representation
$$
D_m[\beta]=\frac{(2m-1)!}{h^{2m}}\Delta_2^{[m]}[\beta]*\sum_{k=1}^{m-1}
\frac{\lambda_k^{|\beta|+m-2}}{E_{2m-2}'(\lambda_k)},
$$
where
$\Delta_2^{[m]}[\beta]={\displaystyle\sum\limits_{k=-m}^m(-1)^{k+m}{2m\choose
m+k}\delta[\beta-k] }$ is symmetric difference of order $2m$.}

{\bf Property 2.} {\it The operator  $D_m[\beta]$ and  monomials
$[\beta]^k=(h\beta)^k$ are connected as}
$$
\sum_{\beta}D_m[\beta][\beta]^k= \left\{
\begin{array}{lll}
0&\mbox{ for } & 0\leq k\leq 2m-1,\\
(2m)!&\mbox{ for } & k= 2m,
\end{array}
\right. \eqno (2)
$$
$$
\sum_{\beta}D_m[\beta][\beta]^k= \left\{
\begin{array}{lll}
0&\mbox{ for } & 2m+1\leq k\leq 4m-1,\\
{\displaystyle             \frac{h^{2m}(4m)!B_{2m}}{(2m)!}    }
&\mbox{ for }  & k= 4m.
\end{array}
\right. \eqno (3)
$$

{\bf Property 3.} {\it The operator $D_m[\beta]$ and the function
$\exp(2\pi i hp\beta)$ connected as
$$
\sum_{\beta}D_m[\beta]\exp(2\pi ihp\beta)=
$$
$$
\frac{(-1)^m2^{2m}(2m-1)!h^{-2m}\sin^{2m}(\pi hp)}
{2\sum\limits_{k=0}^{m-2}a_k^{(2m-2)}\cos 2\pi
hp(m-1-k)+a_{m-1}^{(2m-2)}},
$$
where
$$
a_k^{(2m-2)}=\sum\limits_{j=0}^k(-1)^j{2m\choose j} (k+1-j)^{2m-1}
$$
are the coefficients of the Euler polynomial $E_{2m-2}(\lambda)$.}

\section{Lemmas.}

As known, Euler polynomials $E_k(\lambda)$ have following form
$$
\lambda
E_k(\lambda)=(1-\lambda)^{k+2}D^k\frac{\lambda}{(1-\lambda)^2},
\eqno (4)
$$
where
$$
D=\lambda\frac{d}{d\lambda},\ \
D^k=\lambda\frac{d}{d\lambda}D^{k-1}.
$$
In [8] was shown, that all roots $\lambda_j^{(k)}$ of the Euler
polynomial $E_k(\lambda)$ are real, negative and different:
$$
\lambda_1^{(k)}<\lambda_2^{(k)}<...<\lambda_k^{(k)}<0. \eqno (5)
$$
Furthermore, the roots, equal standing  from the ends of the chain
(5) mutually inverse:
$$
\lambda_j^{(k)}\cdot \lambda_{k+1-j}^{(k)}=1. \eqno (6)
$$
If we denote $E_k(\lambda)=\sum\limits_{s=0}^ka_s^{(k)}\lambda^s$,
then the coefficients  $a_s^{(k)}$ of Euler polynomials, as this
was shown by Euler himself, are expressed by formula
$$
a_s^{(k)}=\sum\limits_{j=0}^s(-1)^j{k+2\choose j} (s+1-j)^{k+1}.
$$
From the definition $E_k(\lambda)$ follow following statement.

{\bf Lemma 1.} \textit{For polynomial $E_k(\lambda)$ following
recurrence relation is valid}
$$
E_k(\lambda)=(k\lambda+1)E_{k-1}(\lambda)+\lambda(1-\lambda)E_{k-1}'(\lambda),
\eqno (7)
$$
\emph{where} $E_0(\lambda)=1$, $k=1,2,....$

{\bf Lemma 2.} \emph{The polynomial  $E_k(\lambda)$ satisfies the
identity}
$$
E_k(\lambda)=\lambda^kE_k\left(\frac{1}{\lambda}\right), \eqno (8)
$$
\emph{or otherwise} $a_s^{(k)}=a_{k-s}^{(k)},\ \ s=0,1,2,...,k$.

\textbf{Proof of lemma 1.} From (4) we can see, that
$$
E_{k-1}(\lambda)=\lambda^{-1}(1-\lambda)^{k+1}D^{k-1}\frac{\lambda}{(1-\lambda)^2}.
\eqno (9)
$$
Differentiating by $\lambda$ the polynomial $E_{k-1}(\lambda)$, we
get
$$
E_{k-1}'(\lambda)=-(1-\lambda)^k\lambda^{-2}(k\lambda+1)D^{k-1}\frac{\lambda}
{(1-\lambda)^2}+\frac{E_k(\lambda)}{\lambda(1-\lambda)}.
$$
Hence and from (9) we obtain, that
$$
(k\lambda+1)E_{k-1}(\lambda)+\lambda(1-\lambda)E_{k-1}'(\lambda)=(k\lambda+1)
\lambda^{-1}(1-\lambda)^{k+1}D^{k-1}\frac{\lambda}{(1-\lambda)^2}-
$$
$$
-(1-\lambda)^{k+1}\lambda^{-1}(k\lambda+1)D^{k-1}\frac{\lambda}{(1-\lambda)^2}+
E_k(\lambda)=E_k(\lambda).
$$
So, lemma 1 is proved.

\textbf{Proof of lemma 2.} The lemma we will prove by induction
method. When $k=1$ from (4) we find
$$
E_1(\lambda)=\lambda+1.
$$
We suppose, that when  $k\geq 1$ the equality
$a_n^{(k-1)}=a_{k-1-n}^{(k-1)}$,
 $n=0,1,...,k-1$ is fulfilled. We assume, that $a_n^{(k-1)}=0$ for  $n<0$ and $n>k-1$.

From (7) we have
 $$
 a_s^{(k)}=(s+1)a_s^{(k-1)}+(k-s+1)a_{s-1}^{(k-1)},
 $$
then, using assumption of induction, we get
 $$
 a_{k-s}^{(k)}=(k-s+1)a_{k-s}^{(k-1)}+(s+1)a_{k-s-1}^{(k-1)}=
 (k-s+1)a_{s-1}^{(k-1)}+(s+1)a_s^{(k-1)}=a_s^{(k)},
 $$
 and lemma 2 is proved.

\section{Proof of theorem 1}

For this we will use function
$$
G_m(x)={x^{2m-1}sign x \over 2\cdot (2m-1)!}.
$$
To this function we correspond following function of discrete
argument:
$$
G_m[\beta]={(h\beta)^{2m-1}sign(h\beta)\over 2\cdot (2m-1)!}.
$$
Here we must find such function $D_m[\beta]$, which satisfies the
equality
$$
h D_m[\beta]*G_m[\beta]=\delta[\beta]. \eqno (10)
$$

According to the theory of periodic generalized functions and
Fourier transformation in them instead of function  $D_m[\beta]$
it is convenient to search harrow shaped function [1]
$$
\stackrel{\abd}{D}_m(x)=\sum\limits_{\beta}D_m[\beta]\delta(x-h\beta).
$$

The equality (10) in the class of harrow shaped functions goes to
equation
$$
h\stackrel{\abd}{D}_m(x)*\stackrel{\abd}{G}_m(x)=\delta(x), \eqno
(11)
$$
where
$$
\stackrel{\abd}{G}_m(x)=\sum\limits_{\beta}G_m[\beta]\delta(x-h\beta).
$$

It is known [1], that the class of harrow shaped functions and the
class of functions of discrete variables are isomorphic. So
instead of function of discrete argument $D_m[\beta]$ it is
sufficiently to investigate the function
$\stackrel{\abd}{D}_m(x)$, defining from equation (11).

Later on we need following well known formulas of Fourier
transformation:
$$
F[f(p)]=\int f(x)\exp(2\pi ipx)dx,
$$
$$
F^{-1}[f(p)]=\int f(x)\exp(-2\pi ipx)dx,
$$
$$
F[f(x)*\varphi(x)]=F[f(x)]\cdot F[\varphi(x)],
$$
$$
F[\delta(x)]=1.
$$

Applying to both parts  of (11) Fourier transformation, we get
$$
F[\stackrel{\abd}{D}_m(x)]\cdot F[h\stackrel{\abd}{G}_m(x)]=1.
\eqno (12)
$$
Fourier transform of $h\stackrel{\abd}{G}_m(p)$ is well known
periodic function, given in $R$ with period  $h^{-1}$
$$
F[h\stackrel{\abd}{G}_m(x)]={(-1)^m\over
(2\pi)^{2m}}\sum\limits_{\beta} \frac{1}{|p-h^{-1}\beta|^{2m}},\ \
p\neq h^{-1}\beta. \eqno (13)
$$
This formula  is obtained from the equalities
$$
F[G_m(p)]={(-1)^m\over (2\pi)^{2m}}\frac{1}{|p|^{2m}}\ \ (\mbox{
[1, p. 729]})
$$
and
$$
\stackrel{\abd}{G}_m(x)=G_m(x)\sum_{\beta}\delta(x-h\beta).
$$
Hence, taking into account  (12), we get
$$
F[\stackrel{\abd}{D}_m(p)]=\left[{(-1)^m\over
(2\pi)^{2m}}\sum_{\beta}\frac{1}
{|p-h^{-1}\beta|^{2m}}\right]^{-1}. \eqno (14)
$$

The main properties of this function in multidimensional case,
appearing in construction of discrete analogue of the polyharmonic
operator, were investigated in [1].

We give some of them, which we will use later on.

\textbf{1.} \emph{Zeros of the function
$F[\stackrel{\abd}{D}_m(p)]$ are the points $p=h^{-1}\beta$.}

\textbf{2.} \emph{The function  $F[\stackrel{\abd}{D}_m(p)]$ is
periodic with period $h^{-1}$,  real and analytic for all real
$p$.}

The function $F[\stackrel{\abd}{D}_m(p)]$ can be represented in
the form of Fourier series
$$
F[\stackrel{\abd}{D}_m(p)]=\sum\limits_{\beta}\hat{D}_m[\beta]\exp(2\pi
ih\beta p), \eqno (15)
$$
where
$$
\hat{D}_m[\beta]=\int\limits_0^{h^{-1}}F[\stackrel{\abd}{D}_m(p)]\exp(-2\pi
ih\beta p)dp. \eqno (16)
$$
Applying inverse Fourier transformation to the equality (15), we
get harrow shaped function
$$
\stackrel{\abd}{D}_m(x)=\sum\limits_{\beta}\hat{D}_m[\beta]\delta(x-h\beta).
\eqno (17)
$$

Thus, $\hat{D}_m[\beta]$ is searching function $D_m[\beta]$ of
discrete argument or discrete analogue of the operator
$\displaystyle\frac{d^{2m}}{dx^{2m}}$. For finding the function
$\hat{D}_m[\beta]$ calculation of the integral (16) inadvisable.
We will find it by following way.

By virtue of known formula
$$
\sum\limits_{\beta}\frac{1}{(p-\beta)^2}=\frac{\pi^2}{\sin^2\pi p}
$$
and from the formula (13) we get
$$
F[h\stackrel{\abd}{G}_1(p)]=\frac{-1}{(2\pi)^2}
\sum\limits_{\beta}\frac{1}{(p-h^{-1}\beta)^2}
=\frac{-h^2}{4\sin^2\pi ph}.
$$
Hence by differentiating we have
$$
\frac{d}{dp}F[h\stackrel{\abd}{G}_1(p)]={2\over (2\pi)^2}
\sum\limits_{\beta}\frac{1}{(p-h^{-1}\beta)^3}.
$$
Thus continuing further, we obtain
$$
\frac{d^{2m-2}}{dp^{2m-2}}F[h\stackrel{\abd}{G}_1(p)]=-{(2m-1)!
\over (2\pi)^{2m}}
\sum\limits_{\beta}\frac{1}{(p-h^{-1}\beta)^{2m}}=
$$
$$
=(-1)^{m-1}(2m-1)!(2\pi)^{2m-2}F[h\stackrel{\abd}{G}_m(p)].
$$
So,
$$
F[h\stackrel{\abd}{G}_m(p)]=\frac{(-1)^mh^2}{2^{2m}\pi^{2m-2}(2m-1)!}
\frac{d^{2m-2}}{dp^{2m-2}}\left(\frac{1}{\sin^2\pi hp}\right).
$$
Consider, the expression
$$
\frac{d^{2m-2}}{dp^{2m-2}}\left(\frac{1}{\sin^2\pi hp}\right).
$$
Using
$$
\sin\pi hp=\frac{\exp(\pi ihp)-\exp(-\pi ihp)}{2i},
$$
we have
$$
\frac{d^{2m-2}}{dp^{2m-2}}\left(\frac{-4}{(\exp(\pi ihp)-\exp(-\pi
ihp))^2}\right)= -4\frac{d^{2m-2}}{dp^{2m-2}}\left(\frac{\exp(2\pi
ihp)}{(\exp(2\pi ihp)-1)^2}\right).
$$
We will do change of variables $\lambda=\exp(2\pi ihp)$, then in
view of that
$$
\frac{d}{dp}=\frac{d\lambda}{dp}\frac{d}{d\lambda}\mbox{ and }
\frac{d}{dp}= 2\pi ih\lambda\frac{d}{d\lambda},
$$
we get
$$
\frac{d^{2m-2}}{dp^{2m-2}}=(2\pi ih)^{2m-2}D^{2m-2},
$$
where
$$
D=\lambda\frac{d}{d\lambda},\ \
D^{2m-2}=\lambda\frac{d}{d\lambda}D^{2m-3}.
$$
Thus,
$$
F[h\stackrel{\abd}{G}_m(p)]=\frac{h^{2m}}{(2m-1)!}D^{2m-2}\frac{\lambda}{(1-\lambda)^2}.
$$
Hence in virtue of (4) we have
$$
F[h\stackrel{\abd}{G}_m(p)]=\frac{h^{2m}}{(2m-1)!}\frac{\lambda
E_{2m-2}(\lambda)} {(1-\lambda)^{2m}}. \eqno (18)
$$
From  (18), according to (12), we obtain
$$
F[\stackrel{\abd}{D}_m(p)]=\frac{(2m-1)!}{h^{2m}}
\frac{(1-\lambda)^{2m}}{\lambda E_{2m-2}(\lambda)}. \eqno (19)
$$
Now in order to obtain Fourier-series expansion, we will do
following.

We divide the polynomial $(1-\lambda)^{2m}$ to the polynomial
$\lambda E_{2m-2}(\lambda)$:
$$
\frac{(1-\lambda)^{2m}}{\lambda\sum\limits_{s=0}^{2m-2}a_s^{(2m-2)}\lambda^s}=
\lambda-2m-a_{2m-3}^{(2m-2)}+ \frac{P_{2m-2}(\lambda)}{\lambda
E_{2m-2}(\lambda)}, \eqno (20)
$$
where $P_{2m-2}(\lambda)$ is a polynomial of degree $2m-2$. It is
not difficult to see, that the rational fraction
$\displaystyle\frac{P_{2m-2}(\lambda)}{\lambda E_{2m-2}(\lambda)}$
is proper fraction, i.e. degree of the polynomial
$P_{2m-2}(\lambda)$ is less than degree of the polynomial $\lambda
E_{2m-2}(\lambda)$. Since the roots of the polynomial
$E_{2m-2}(\lambda)$ are real and different, then the rational
fraction $\displaystyle\frac{P_{2m-2}(\lambda)}{\lambda
E_{2m-2}(\lambda)}$ is expanded to the sum of elementary
fractions. Searching expansion has following form
$$
\frac{P_{2m-2}(\lambda)}{\lambda E_{2m-2}(\lambda)}=
{A_0\over\lambda}+
\sum\limits_{k=1}^{m-1}\frac{A_{1,k}}{\lambda-\lambda_{1,k}}+
\sum\limits_{k=1}^{m-1}\frac{A_{2,k}}{\lambda-\lambda_{2,k}},
\eqno (21)
$$
where $A_0,\ \ A_{1,k}, \ \ A_{2,k}$ are unknown coefficients,
$\lambda_{1,k}$ are the roots of the polynomial
$E_{2m-2}(\lambda)$, in modulus less than unity, and
$\lambda_{2,k}$ are the roots of the polynomial
$E_{2m-2}(\lambda)$, in modulus greater than unity. By (21) the
equality (20) takes the form
$$
\frac{(1-\lambda)^{2m}}{\lambda\sum\limits_{s=0}^{2m-2}a_s^{(2m-2)}\lambda^s}=
\lambda-2m-a_{2m-3}^{(2m-2)}+ {A_0\over\lambda}+
$$
$$
+\sum\limits_{k=1}^{m-1}\frac{A_{1,k}}{\lambda-\lambda_{1,k}}+
\sum\limits_{k=1}^{m-1}\frac{A_{2,k}}{\lambda-\lambda_{2,k}}.
\eqno (22)
$$
Reducing to the common denominator and omitting it, we get
$$
(1-\lambda)^{2m}=\lambda^2E_{2m-2}(\lambda)-\lambda(2m+a_{2m-3}^{(2m-2)})E_{2m-2}(\lambda)+
$$
$$
+A_0E_{2m-2}(\lambda)+ \sum\limits_{k=1}^{m-1}\frac{A_{1,k}\lambda
E_{2m-2}(\lambda)}{\lambda-\lambda_{1,k}}+
\sum\limits_{k=1}^{m-1}\frac{A_{2,k}\lambda
E_{2m-2}(\lambda)}{\lambda-\lambda_{2,k}}. \eqno (23)
$$
Assuming in the equality  (23) consequently $\lambda=0$,
$\lambda=\lambda_{1,k}$ and $\lambda=\lambda_{2,k}$, we find
$$
1=E_{2m-2}(0)A_0;\ \
(1-\lambda_{1,k})^{2m}=\lambda_{1,k}E_{2m-2}'(\lambda_{1,k})A_{1,k};
$$
$$
(1-\lambda_{2,k})^{2m}=\lambda_{2,k}E_{2m-2}'(\lambda_{2,k})A_{2,k}.
$$
Hence
$$
A_0=1;\ \
A_{1,k}=\frac{(1-\lambda_{1,k})^{2m}}{\lambda_{1,k}E_{2m-2}'(\lambda_{1,k})};\
\
A_{2,k}=\frac{(1-\lambda_{2,k})^{2m}}{\lambda_{2,k}E_{2m-2}'(\lambda_{2,k})}.
$$
Using (6), we have
$$
A_{1,k}=\frac{(1-{1\over
\lambda_{2,k}})^{2m}}{\lambda_{2,k}^{-1}E_{2m-2}'({1\over
\lambda_{2,k}})}=
\frac{(\lambda_{2,k}-1)^{2m}}{\lambda_{2,k}^{2m-1}E_{2m-2}'({1\over
\lambda_{2,k}})}.
$$
In virtue of (7) we obtain
$$
{1\over \lambda_{2,k}}\left(1-{1\over \lambda_{2,k}}\right)
E_{2m-2}'({1\over \lambda_{2,k}})=E_{2m-1}({1\over
\lambda_{2,k}}),
$$
$$
\lambda_{2,k}(1-\lambda_{2,k})E_{2m-2}'(\lambda_{2,k})=E_{2m-1}(\lambda_{2,k}),
$$
hence
$$
E_{2m-2}'({1\over \lambda_{2,k}})=\frac{E_{2m-1}({1\over
\lambda_{2,k}}) \lambda_{2,k}^2}{\lambda_{2,k}-1},
$$
$$
E_{2m-2}'(\lambda_{2,k})=\frac{E_{2m-1}(\lambda_{2,k})}{\lambda_{2,k}
(1-\lambda_{2,k})}.
$$
From here application of the lemma 2 gives
$$
A_{1,k}={-A_{2,k}\over\lambda_{2,k}^2},\ \
A_{1,k}={(1-\lambda_{1,k})^{2m+1}\over E_{2m-1}(\lambda_{1,k})}.
\eqno (24)
$$
Since $|\lambda_{1,k}|<1$ and $|\lambda_{2,k}|>1$, then
$$
\sum\limits_{k=1}^{m-1}\frac{A_{1,k}}{\lambda-\lambda_{1,k}}
\mbox{ and }
\sum\limits_{k=1}^{m-1}\frac{A_{2,k}}{\lambda-\lambda_{2,k}}
$$
can be represented as Laurent series on the circle
$|\lambda^2|=1$:
$$
\sum\limits_{k=1}^{m-1}\frac{A_{1,k}}{\lambda-\lambda_{1,k}}=
{1\over \lambda}\sum\limits_{k=1}^{m-1}\frac{A_{1,k}}{1-
{\lambda_{1,k}\over \lambda}}= {1\over
\lambda}\sum\limits_{k=1}^{m-1}A_{1,k}
\sum_{\beta=0}^{\infty}\left(\lambda_{1,k}\over
\lambda\right)^{\beta}, \eqno (25)
$$
$$
\sum\limits_{k=1}^{m-1}\frac{A_{2,k}}{\lambda-\lambda_{2,k}}=
-\sum\limits_{k=1}^{m-1}\frac{A_{2,k}}{\lambda_{2,k}(1-
{\lambda\over \lambda_{2,k}})}=
-\sum\limits_{k=1}^{m-1}\frac{A_{2,k}}{\lambda_{2,k}}
\sum_{\beta=0}^{\infty}\left(\lambda\over
\lambda_{2,k}\right)^{\beta}. \eqno (26)
$$
Putting (25), (26) to (22) and taking into account
  $\lambda=\exp(2\pi ihp)$ from (19), (20), we obtain
$$
F[\stackrel{\abd}{D}_m(p)]=\frac{(2m-1)!}{h^{2m}}\Bigg[\exp(2\pi
ihp)-2m- a_{2m-3}^{(2m-2)}+
$$
$$
+\exp(-2\pi ihp)+\sum\limits_{k=1}^{m-1}\Bigg( A_{1,k}
\sum_{\beta=0}^{\infty}\lambda_{1,k}^{\beta}\exp(-2\pi
ihp(\beta+1))-
$$
$$
-A_{2,k}\lambda_{2,k}^{-1}\sum_{\beta=0}^{\infty}
\left({1\over\lambda_{2,k}}\right)^{\beta}\exp(2\pi
ihp\beta)\Bigg)\Bigg].
$$
Thus, searching Fourier series for $F[\stackrel{\abd}{D}_m(p)]$
have following form
$$
F[\stackrel{\abd}{D}_m(p)]=\sum_{\beta}D_m[\beta]\exp(2\pi ih\beta
p),
$$
where
$$
D_m[\beta]=\frac{(2m-1)!}{h^{2m}}\left\{
\begin{array}{lll}
\sum\limits_{k=1}^{m-1}A_{1,k}\lambda_{1,k}^{-\beta-1}& \mbox{ for }& \beta\leq -2,\\
1+\sum\limits_{k=1}^{m-1}A_{1,k}& \mbox{ for }& \beta= -1,\\
-2^{2m-1}-\sum\limits_{k=1}^{m-1}A_{2,k}\lambda_{2,k}^{-1}& \mbox{ for }& \beta= 0,\\
1-\sum\limits_{k=1}^{m-1}A_{2,k}\lambda_{2,k}^{-2}& \mbox{ for }& \beta= 1,\\
-\sum\limits_{k=1}^{m-1}A_{2,k}\lambda_{2,k}^{-\beta-1}& \mbox{
for }& \beta\geq 2.
\end{array}
\right.
$$
With the help (24) the function $D_m[\beta]$ we rewrite in the
form
$$
D_m[\beta]=\frac{(2m-1)!}{h^{2m}}\left\{
\begin{array}{lll}
{\displaystyle
\sum\limits_{k=1}^{m-1}\frac{(1-\lambda_{1,k})^{2m+1}\lambda_{1,k}^{|\beta|}}
{\lambda_kE_{2m-1}(\lambda_{1,k})}      }
& \mbox{ for }& |\beta|\geq 2,\\
{\displaystyle
1+\sum\limits_{k=1}^{m-1}\frac{(1-\lambda_{1,k})^{2m+1}}
{E_{2m-1}(\lambda_{1,k})}        }
& \mbox{ for }& |\beta|= 1,\\
{\displaystyle-2^{2m-1}+\sum\limits_{k=1}^{m-1}\frac{(1-\lambda_{1,k})^{2m+1}}
{\lambda_{1,k}E_{2m-1}(\lambda_{1,k})}          } & \mbox{ for }&
\beta= 0.
\end{array}
\right.
$$
We note, that
$$
D_m[\beta]=D_m[-\beta].
$$

Theorem  1 is proved completely.

\section{Proofs of properties.}

\textbf{Proof of property 1.}

Following is takes placed
$$
F[\stackrel{\abd}{\Delta}^{[1]}_2(p)]=-4\sin^2\pi ph.
$$
Indeed,
$$
\stackrel{\abd}{\Delta}^{[1]}_2(x)=\stackrel{\abd}{\delta}(x+1)-
2\stackrel{\abd}{\delta}(x)+2\stackrel{\abd}{\delta}(x-1)=\delta(x+h)-
2\delta(x)+\delta(x-h).
$$
By definition of Fourier transformation we have
$$
F[\stackrel{\abd}{\Delta}^{[1]}_2(p)]=\int\exp(2\pi ipx)
\stackrel{\abd}{\Delta}^{[1]}_2(x)dx=-4\sin^2\pi ph.
$$

Hence consequently we obtain
$$
F[\stackrel{\abd}{\Delta}^{[m]}_2(p)]=
F[\overbrace{\stackrel{\abd}{\Delta}^{[1]}_2(p)*\stackrel{\abd}{\Delta}^{[1]}_2(p)*...*
\stackrel{\abd}{\Delta}^{[1]}_2(p)}^{m \mbox{ times}}]=
(-4)^m\sin^{2m}\pi ph. \eqno (27)
$$
Immediately we have
$$
\left[\frac{(1-\lambda)^2}{-4\lambda}\right]^m=\sin^{2m}\pi hp.
\eqno (28)
$$
By virtue of (27) and (28) the formula (19) takes form
$$
F[\stackrel{\abd}{D}_m(p)]=\frac{(2m-1)!}{h^{2m}}\frac{\lambda^{m-1}}{E_{2m-2}(\lambda)}
F[\stackrel{\abd}{\Delta}^{[m]}_2(p)]. \eqno (29)
$$
Now expanding rational fraction ${\lambda^{m-1}\over
E_{2m-2}(\lambda)}$ to the sum of elementary fractions, we have
$$
{\lambda^{m-1}\over
E_{2m-2}(\lambda)}=\sum\limits_{k=1}^{m-1}\left[
\frac{B_{1,k}}{\lambda-\lambda_{1,k}}+
\frac{B_{2,k}}{\lambda-\lambda_{2,k}}\right], \eqno (30)
$$
where
$$
B_{1,k}=\frac{\lambda_{1,k}^{m-1}}{E_{2m-2}'(\lambda_{1,k})},\ \
B_{2,k}=\frac{\lambda_{2,k}^{m-1}}{E_{2m-2}'(\lambda_{2,k})}.
$$
Since $|\lambda_{1,k}|<1$ and $|\lambda_{2,k}|>1$, then expanding
$\frac{B_{1,k}}{\lambda-\lambda_{1,k}}$ and
$\frac{B_{2,k}}{\lambda-\lambda_{2,k}}$ to the Laurent series on
the circle $|\lambda|=1$, we find
$$
\frac{B_{1,k}}{\lambda-\lambda_{1,k}}={1\over\lambda}\cdot
\frac{B_{1,k}}{1-{\lambda_{1,k}\over\lambda}}={B_{1,k}\over\lambda}
\sum\limits_{\beta=0}^{\infty}\left({\lambda_{1,k}\over\lambda}\right)^{\beta},
\eqno (31)
$$
$$
\frac{B_{2,k}}{\lambda-\lambda_{2,k}}=
-\frac{B_{2,k}}{\lambda_{2,k}(1-{\lambda\over\lambda_{2,k}})}=
-{B_{2,k}\over\lambda_{2,k}}\sum\limits_{\beta=0}^{\infty}\left({\lambda\over\lambda_{2,k}}\right)^{\beta}.
\eqno (32)
$$
On the strength of (6), (7), (8), (31), (32) the equality (30)
takes the form
$$
{\lambda^{m-1}\over E_{2m-2}(\lambda)}=
\sum\limits_{k=1}^{m-1}\frac{\lambda_k^{m-2}}{E_{2m-2}'(\lambda_k)}
\sum\limits_{\beta}\lambda_k^{|\beta|}\lambda^{\beta}, \ \
\lambda_k=\lambda_{1,k}.
$$
So, from (29) we get
$$
F[\stackrel{\abd}{D}_m(p)]=\frac{(2m-1)!}{h^{2m}}\sum\limits_{k=1}^{m-1}
\frac{\lambda_k^{m-2}}{E_{2m-2}'(\lambda_k)}
\sum\limits_{\beta}\lambda_k^{|\beta|}\lambda^{\beta}
F[\stackrel{\abd}{\Delta}_2^{[m]}(p)]. \eqno (33)
$$
Applying to (33) inverse Fourier transformation, after some
simplifications we obtain
$$
\stackrel{\abd}{D}_m(x)=\frac{(2m-1)!}{h^{2m}}\sum\limits_{\beta}
\Delta_2^{[m]}[\beta]*\sum\limits_{k=1}^{m-1}\frac{\lambda_k^{|\beta|+m-2}}
{E_{2m-2}'(\lambda_k)}\delta(x-h\beta)=
$$
$$
=\sum\limits_{\beta}D_m[\beta]\delta(x-h\beta). \eqno (34)
$$
According to definition of harrow shaped functions from (34) we
get the statement of  property 1.

\textbf{Proof of property 2.}

The equality (2) proved in [1].\\
Here we will prove (3).

From (15), (17) we obtain
$$
F[\stackrel{\abd}{D}_m(p)]=\sum\limits_{\beta}D_m[\beta]\exp(2\pi
ih\beta p). \eqno (35)
$$
Using expansion of $\exp(2\pi ihp\beta)$ and by formula (2) in
cases $0\leq k\leq 2m-1$ and $k=2m$, we will have
$$
F[\stackrel{\abd}{D}_m(p)]=\sum\limits_{\beta}D_m[\beta]
\sum_{k=0}^{\infty}\frac{(2\pi ih\beta p)^k}{k!}=
$$
$$
= \sum_{\beta}D_m[\beta]\sum\limits_{k=2m+1}^{\infty}\frac{(2\pi
ihp\beta)^k}{k!}+ (2\pi ip)^{2m}. \eqno (36)
$$
On the other hand, from (14) we get
$$
F[\stackrel{\abd}{D}_m(p)]= (2\pi i)^{2m}\left[\sum\limits_{\beta}
{1\over (p-h^{-1}\beta)^{2m}}\right]^{-1}. \eqno (37)
$$
Thus,  on the strength of (36) and (37) we have
$$
\sum_{\beta}D_m[\beta]\sum\limits_{k=2m+1}^{\infty}\frac{(2\pi
ihp\beta)^k}{k!}+ (2\pi ip)^{2m}=(2\pi
i)^{2m}\left[\sum\limits_{\beta} {1\over
(p-h^{-1}\beta)^{2m}}\right]^{-1},
$$
or
$$
\sum_{\beta}D_m[\beta]\sum\limits_{k=2m+1}^{\infty}\frac{(2\pi
ihp\beta)^k}{k!}= -(2\pi
ip)^{2m}\left(1-\left[\sum\limits_{\gamma} {(ph)^{2m}\over
(ph-\gamma)^{2m}}\right]^{-1}\right). \eqno (38)
$$
Consider
$$
\psi(h,p,m)=1-\left[\sum\limits_{\gamma}\frac{(ph)^{2m}}{(ph-\gamma)^{2m}}\right]^{-1}=
$$
$$
=1-\left[1+(ph)^{2m}\sum_{\gamma=1}^{\infty}\left[\frac{1}{(ph-\gamma)^{2m}}+
\frac{1}{(ph+\gamma)^{2m}}\right]\right]^{-1}=
$$
$$
=1-\left[1+ (ph)^{2m}\sum_{\gamma=1}^{\infty}\gamma^{-2m}
\left[(1-\frac{ph}{\gamma})^{-2m}+
(1+\frac{ph}{\gamma})^{-2m}\right] \right]^{-1}.
$$
Hence choosing  $p$ such, that
$$
Q= (ph)^{2m}\sum_{\gamma=1}^{\infty}\gamma^{-2m}
\left[(1-\frac{ph}{\gamma})^{-2m}+
(1+\frac{ph}{\gamma})^{-2m}\right]<1,
$$
expanding the fraction $\displaystyle\frac{1}{1+Q}$ to the series
of geometric progression, we have
$$
\psi(h,p,m)=(ph)^{2m}\sum_{\gamma=1}^{\infty}\gamma^{-2m}
\left[(1-\frac{ph}{\gamma})^{-2m}+
(1+\frac{ph}{\gamma})^{-2m}\right]+
$$
$$
+O((hp)^{4m})=2(ph)^{2m}\sum_{\gamma=1}^{\infty}\gamma^{-2m}+O(h^{2m+1}).
\eqno (39)
$$
The left part of  (38) we rewrite in the form
$$
\sum\limits_{\beta}h^{2m}D_m[\beta]\sum_{k=2m+1}^{\infty}h^{k-2m}
\frac{(2\pi ip\beta)^k}{k!}=\sum_{n=1}^{\infty}h^n\mu(p,n,m),
\eqno (40)
$$
where
$\mu(p,n,m)={\displaystyle\sum\limits_{\beta}h^{2m}D_m[\beta]\frac{(2\pi
ip\beta)^{2m+n}} {(2m+n)!}}$ does not depend on $h$, since
$h^{2m}D_m[\beta]$ does not depend  $h$, that clear from (1).

Comparing right hand sides of (39) and (40), from (38) we have
$\mu(p,n,m)=0$ for $n=1,2,...,2m-1$, i.e. for $k=2m+1,...,4m-1$
$$
\mu(p,2m,m)=\sum\limits_{\beta}h^{2m}D_m[\beta]\frac{(2\pi
ip\beta)^{4m}}{(4m)!}=
(-1)^{m+1}2(2\pi)^{2m}p^{4m}\sum\limits_{\gamma=1}^{\infty}\gamma^{-2m}.
$$
Hence after some calculations and by virtue of the equality
$$
\sum\limits_{\gamma=1}^{\infty}\gamma^{-2m}=\frac{(-1)^{m-1}(2\pi)^{2m}B_{2m}}
{2\cdot (2m)!}
$$
we get
$$
\sum_{\beta}D_m[\beta][\beta]^{4m}=\frac{h^{2m}(4m)!B_{2m}}{(2m)!},
$$
which proves the property 2.

\textbf{Proof of property 3.} From (19) and (35) we have
 $$
\sum_{\beta}D_m[\beta]\exp(2\pi
ihp\beta)=\frac{(2m-1)!(1-\lambda)^{2m}} {h^{2m}\lambda
E_{2m-2}(\lambda)}.
$$
Using (8), (40) and $\lambda=\exp(2\pi ihp\beta)$, after some
simplifications we get property 3.

\begin{center}
\textbf{REFERENCES}
\end{center}
\begin{enumerate}
\item Sobolev S.L. Introduction to the Theory of Cubature
Formulas. -Moskow.: Nauka, 1974. - 808 p.

\item Zhamalov Z.Zh. About one problem of Wiener-Hopf appearing
in optimization of quadrature formulas. - In the Book: Boundary
problems for differential equations. - Tashkent: Fan, 1975. - pp.
129-150.

\item Zhamalov Z.Zh. About one difference analogue of operator
$\frac{d^{2m}}{dx^{2m}}$ and its construction. - In the Book:
Direct and inverse problem for differential equations with partial
derivatives and their properties. -Tashkent:  Fan, 1978. - pp.
97-108.

\item Shadimetov Kh.M. Optimal quadrature formulas in the
$L_2^{(m)}(\Omega)$ and $L_2^{(m)}(R^1)$ // Dokl. Akad. Nauk.
Uzbekistan. - 1983. - No.3. - pp. 5-8.

\item Shadimetov Kh.M. Discrete analogue of the
$\frac{d^{2m}}{dx^{2m}}$ and its construction // Questions of
Computational and Applied Mathematics. - Tashkent, 1985. - pp.
22-35.

\item Shadimetov Kh.M. About one explicit representation of
discrete analogue of the differential operator of order $2m$
// Dokl. Akad. Nauk. Uzbekistan. - 1996. - No.9. - pp. 5-7.

\item Shadimetov Kh.M. About one method of solution of
difference equation with convolution for computation of optimal
coefficients of Sobolev's quadrature formulas // Uzbek
Mathematical Journal. - 1998. - No.4. - pp. 68-76.

\item Frobenius.  Uber Bernoullische Zahlen und Eulersche
Polynomen Sitzungsberichte der Preus\-sischen Akademic de
Wissenschaften. 1910.
\end{enumerate}

\textbf{Shadimetov Kholmat Mahkambaevich}\\
Head of the Department of Computational Methods\\
Institute of Mathematics and Information Technologies\\
Tashkent, Uzbekistan\\

\end{document}